\documentclass[11pt]{amsart}
\usepackage{amssymb}
\usepackage[english]{babel} 
\def\R{\mathbb R}

\def\N{\mathbb N}

\newtheorem{theoreme}{Theorem}

\newtheorem{remark}{Remark}

\newtheorem{proposition}{Proposition}
\newtheorem{lemme}[proposition]{Lemma}

\newtheorem{remarque}[proposition]{Remark}

\numberwithin{equation}{section}
\numberwithin{proposition}{section}

\begin{document}
\title[On the KP-I equation]
{On finite energy solutions of the KP-I equation}
\author{H. Koch}
\address{Universit\"at Dortmund, 44221 Dortmund}
\email{herbert.koch@mathematik.uni-dortmund.de}
\thanks{The work of the first author was partly supported by MSRI and the 
Miller Institute of Basic Research.}   
\author{N. Tzvetkov}
\address{
D\'epartement de Math\'ematiques, Universit\'e Lille I, 59 655 Villeneuve d'Ascq Cedex
}
\email{nikolay.tzvetkov@math.univ-lille1.fr}
\selectlanguage{english}
\begin{abstract}
{We prove that the flow map of the Kadomtsev-Petviashvili-I (KP-I) equation is
not uniformly continuous on bounded sets of the natural energy space.} 
\end{abstract} 
\subjclass{ 35Q55, 37K10}
\keywords{dispersive equations, nonlinear waves.}
\maketitle 
\section{Introduction}
The understanding of solutions to dispersive equations has considerably
deepened during the last decade. Much of the progress is based on the
idea of Bourgain of using $L^2$ based function spaces adapted to the
linear operator. This technique establishes often  a
 very nontrivial  domination of the nonlinearity by the linear operator, 
at least on small scales. It is connected to objects in harmonic
analysis: restriction theorems, local smoothing, maximal functions and 
multilinear estimates. If applicable it leads to existence and
uniqueness via Picard iteration or the implicit function theorem, and
hence to uniform continuity and even differentiability of the flow map.  

Despite the amazing success of this approach there are several problems 
where it failed completely. Two of the most interesting of them are the
Benjamin-Ono  and the Kadomtsev-Petviashvili-I (KP-I)  equation. The Benjamin-Ono equation has been intensively studied during the last three years, and a very
 precise understanding has emerged: In the same way as for Burgers equation the low frequency part leads to a change of the speed of waves (\cite{KT}), 
which in turn contradicts uniform continuity of the flow map. This, however, does not imply illposedness. The transport  effect can be 
controlled by a gauge transform, see Tao \cite{tao}. Ionescu and Kenig \cite{IK} approached 
the transformed problem   by adapted function spaces and bilinear estimates and obtained well-posedness for initial data in $L^2$. 

In this paper we study the KP-I equation
\begin{equation}\label{1}
u_{t}+u_{xxx}-\partial_{x}^{-1}u_{yy}+uu_{x}=0,
\end{equation}
where $(t,x,y)\in\R^3$, $u$ is a real valued function and $\partial_{x}^{-1}$
is a formal notation for the antiderivative, which always exists for tempered distributions, but whose uniqueness requires further considerations. In this paper we only deal with antiderivatives in $L^2$ of $L^2$ functions $f$. In this case 
the antiderivative is uniquely defined and  its Fourier transform can be defined through  multiplication of $\widehat{f}$ by $(i \xi)^{-1}$. It is not hard to see that, if $f$ is in addition compactly supported,  then it  is integrable 
with mean zero, and hence the antiderivative could also be defined through the indefinite integral from $-\infty$.

The KP-I equation appears as an asymptotic model for the propagation of
long, essentially one directional, small amplitude
surface waves when the surface tension is bigger than some critical value.
For smaller values of the surface tension we get the KP-II model.
The KP-I equation can be written in the Lax pair form (see \cite{ZS}) and
thus it shares many features with the "integrable PDE's". 
One also has a family of particular solitary waves solutions called lump solutions. The
study of the KP-I flow close to the lumps is a challenging issue.

A simple calculation shows  that the solutions of (\ref{1}) satisfy, at least
formally, the conservation of the $L^2$ norm
$$
N(u)=\int_{\R^2}u^{2}(t,x,y)dxdy={\rm const}
$$
and the conservation of the energy
$$
E(u)=\frac{1}{2}\int_{\R^2}\Big[(\partial_{x}u)^{2}+(\partial_{x}^{-1}u_{y})^{2}-\frac{1}{3}u^{3}\Big](t,x,y)dxdy={\rm
const}\,.
$$
 Taking into account the anisotropic
Sobolev inequality (see e.g. \cite{Tom} and Lemma \ref{sobolev}),
$$
\|u\|_{L^3(\R^2)}^{3}\leq C
\|u\|_{L^2(\R^2)}^{\frac{3}{2}}\|u_x\|_{L^2(\R^2)}\|\partial_{x}^{-1}u_{y}\|_{L^2(\R^2)}^{\frac{1}{2}}
$$
we deduce that the subspace of $L^2$ of finite energy
is a  Banach space $X$, which we call {\em energy space}, equipped with the norm
$$
\|u\|_{X}=\|u\|_{L^2(\R^2)}+\|u_x\|_{L^2(\R^2)}+\|\partial_{x}^{-1}u_y\|_{L^2(\R^2)}\,.
$$
This provides a natural framework to study the nonlinear problem (\ref{1}). 
The Cauchy problem for the KP-I equation is known to be globally well-posed in
spaces smaller than the energy space (see \cite{MST1,K}). More precisely, the
Cauchy problem associated to (\ref{1}) is globally well-posed for data in the
space $Z$ equipped with the norm
$$
\|u\|_{Z}=\|u\|_{L^2(\R^2)}+\|u_{xx}\|_{L^2(\R^2)}+\|\partial_{x}^{-2}u_{yy}\|_{L^2(\R^2)}\,.
$$
Notice that $u\in Z$ implies $u_x\in L^2(\R^2)$, $u_y\in L^2(\R^2)$,
$\partial_{x}^{-1}u_{y}\in L^2(\R^2)$. Moreover, 
the Lax pair formulation of the KP-I equation implies that
if the initial data
$$
u_0=u|_{t=0}\in Z
$$ 
of (\ref{1}) is smooth (i.e. in the intersection of all
$H^s(\R^2)$, $s\in \N$), and if $\partial_{x}^{-1}u_0$ is also smooth then the
global solution of (\ref{1}) with data $u_0$ satisfies a third conservation
law.
Namely
$$
F(u)={\rm const},
$$
where
\begin{equation}\label{F}
\begin{split}
F(u)  = &\frac{3}{2} \int_{\R^2} u_{xx}^2 + 5 \int_{\R^2} u_y^2
+ \frac{5}{6}\int_{\R^2} (\partial^{-2}_x u_{yy} )^2
\\
 & 
-\frac{5}{6} \int_{\R^2} u^2 (\partial^{-2}_x u_{yy})
-\frac{5}{6} \int_{\R^2} u\, (\partial_{x}^{-1} u_y)^2
 +\frac{5}{4} \int_{\R^2} u^2\, u_{xx} + \frac{5}{24} \int_{\R^2} u^4 \,.
\end{split}
\end{equation}
There is in fact an infinite sequence of formal conservation laws associated to the
KP-I equation (see \cite{ZS}). However, as noticed in \cite{MST1}, it is hard
to find a suitable framework of distributions on $\R^2$ where these
conservation laws make sense.

\medskip

It is presently not known whether  (\ref{1}) is well-posed in the energy space
$X$, but we hope this problem will be given an affirmative answer in the near
future. The goal of this paper
is to show that, whatever the answer is, the flow map of (\ref{1}) can not
be uniformly continuous on bounded sets of the energy space $X$.
Recall that, if one solves (\ref{1}) in $X$ by the Picard iteration, then the
flow map is automatically uniformly continuous on bounded subsets of $X$. 
Our result thus implies  that the solution to (\ref{1}) cannot 
by constructed by  Picard iteration scheme, in sharp contrast with many other
dispersive models as the KdV equation, the KP-II equation (where
(\ref{1})  $-\partial_{x}^{-1}u_{yy}$ is replaced by
$\partial_{x}^{-1}u_{yy}$), the  KdV equation etc...
This feature of the KP-I equation was already observed in \cite{MST1,MST2}.
In the present paper we construct some solutions of (\ref{1}) which are
``responsible'' for this phenomenon.
Here is the precise statement of our result.
\begin{theoreme}\label{th1}
There exist two positive constants $c$ and $C$ and 
two sequences $(u_n)$ and $(\widetilde{u}_{n})$ of solutions of (\ref{1})  
such that for every $t\in [-1,1]$,
$$
\sup_{n}\|u_{n}(t,\cdot)\|_{X}+\sup_{n}\|\widetilde{u}_{n}(t,\cdot)\|_{X}\leq C\, ,
$$
$(u_n)$ and $(\widetilde{u}_{n})$ satisfy initially
$$
\lim_{n\rightarrow\infty}\|u_{n}(0,\cdot)-\widetilde{u}_{n}(0,\cdot)\|_{X}=0,
$$
but, for every $t\in [-1,1]$,
\begin{equation*}
\liminf_{n\rightarrow \infty}\|u_{n}(t,\cdot)-\widetilde{u}_{n}(t,\cdot)\|_{X}\geq
c\,|t|\, .
\end{equation*}
\end{theoreme}
In a previous paper \cite{KT}, we proved a similar result for the Benjamin-Ono
equation. The analysis in the KP-I context is more involved since we use an
additional cancellation in the construction of the approximate solutions,
related to the existence of zero speed waves in the $x$ direction for the
linear KP-I equation. Of course, similar waves do not exist in the KP-II context. In
addition, our analysis uses the Burgers type cancellation which was the only
cancellation involved in the construction of \cite{KT}. Let us also notice
that a technical modification of the proof of Theorem~\ref{th1} is likely to show
that in 
Theorem~\ref{th1} one can replace the energy space by the Sobolev spaces
$H^{s}(\R^2)$, $s>0$ or the spaces $Y_s$ considered in \cite{K} equipped with
the norm
$$
\|u\|_{Y_{s}}=\|u\|_{L^2(\R^2)}+
\|D_x^s u\|_{L^{2}(\R^2)}+
\|\partial_{x}^{-1}u_y\|_{L^2(\R^2)}\,.
$$
\section{Outline of the proof of Theorem~\ref{th1}}
We decompose the proof into three parts: 
\begin{enumerate} 
\item We construct a family of approximate solutions $u_{ap}$ 
depending on parameters $\omega$, 
$|\omega| \le 1$ and $\lambda\gg  1$. Changing $\omega$ leads to 
a phase shift in the high frequency part for positive $t$, but for $t=0$ the
variation of $\omega$ is uniformly smooth. We show that  the residual terms are
small uniformly in all parameters, see Lemma \ref{L2}.
\item We study the bounds of the solutions $u$ with initial datum $u_{ap}(0)$ 
for $t\le 1$ in many  $L^2$ based  spaces. This part relies on the well-posedness for 
smooth data as well as on  conserved quantities.
\item Energy arguments control $\Vert  u_{ap}-u \Vert_{L^2}$. Interpolation 
with the bounds for $u$ and $ u_{ap}$ yields that $u$ is close to $u_{ap}$ 
in suitable function spaces. 
\end{enumerate}
This yields the desired conclusion because $u_{ap}$ depends in a
transparent way on $\omega$, which contradicts  uniform continuity.

Let us explain the idea of the construction of the approximate solution. 
We denote $i$ times the symbol of the spatial part of the linear equation by 
\[ p(\xi,\eta) =   \xi^3 +  \xi^{-1} \eta^2  \] 
Let $\lambda \gg 1$ be a large parameter. The function 
\begin{equation}\label{planewave} 
\cos( \lambda x + 4 \lambda^3 t + \sqrt{3} \lambda^2 y)   
\end{equation}
is a solution to the linear equation. Its velocity vector is 
\begin{equation} \label{velo} 
 \nabla p(\pm \lambda, \pm \sqrt{3} \lambda^2) = \left(\begin{matrix} 
0 \\ \pm 2 \sqrt{3} \lambda \end{matrix}\right). 
\end{equation} 
In particular the velocity of the plane wave \eqref{planewave} 
in the $x$ direction vanishes, which is the reason 
for choosing these points  in the frequency  space.

We fix for the sequel two constants
\begin{equation} \label{alphabeta} 
\frac{1}{2}<\alpha<1<\beta,\quad \alpha+\beta<2\,
\end{equation} 
and  a function $ \varphi \in C_0^{\infty}(\R)$ supported in $[-2,2]$
and identically $1$ in $[-1,1]$.
Since the Fourier transform of 
\[ 
\varphi(x/\lambda^\alpha)\varphi(x/\lambda^\beta)  \cos(  \lambda x + 4 \lambda^3 t
+ \sqrt{3} \lambda^2 y)   
\]
is supported in a small neighborhood of $( \lambda , \sqrt{3} \lambda)$ 
its velocities a close to \eqref{velo}, and hence is an approximate solution. 
Below we will discuss the restrictions on $\alpha$ and $\beta$ needed for this
argument. 

A first guess for the approximate solution is  

\begin{equation} \label{ansatz1}
\begin{split} 
u_{ap}(t,x,y) = &
-\lambda^{-1-\frac{\alpha+\beta}{2}}\varphi(x/\lambda^\alpha)
\varphi(y/\lambda^\beta)\cos(4\lambda^{3}t+\lambda
x+\sqrt{3}\lambda^{2}y + \omega t)
\\
 &
-\lambda^{-1}\omega\, \varphi(x/(2\lambda^\alpha))
\varphi(y/(2\lambda^\beta)).
\end{split} 
\end{equation} 
The crucial point is the dependence of the first part on $\omega$. If we plug 
$u_{ap}$ into the equation then the time derivative leads to an additional term
(compared to $\omega=0$), which is linear in $\omega$, and which essentially 
cancels against the product of the first and the second term in the nonlinearity.

The range of $\alpha$ and $\beta$ is dictated by the following conditions: 
\begin{enumerate} 
\item The low frequency part (the second term in \eqref{ansatz1})
has to converge to $0$ in $L^2$ as $\lambda \to \infty$ uniformly in
      $|\omega|\le 1$.
 Its norm is a constant times 
  \[ |\omega| \lambda^{-1+ {\frac{\alpha+\beta}2}} \]
hence we need $ \alpha + \beta < 2$. Because of the structure of $u_{ap}$ 
convergence of the difference to zero in $X$ follows as well. 
\item The velocity in $y$ direction is of size $\lambda$. Hence $\beta >
      1$ 
is needed so that the high frequency part is confined up to time $1$ in
      an interval of size $\lambda^\beta$ in $y$ direction.
\item Let $\Delta_x= \lambda^\alpha$ be the spatial scale. By the uncertainty
principle the uncertainty in frequency is $\lambda^{-\alpha}$. Then the uncertainty 
in the velocity in $x$ direction is $ \lambda^{1-\alpha}$, leading to 
$\Delta_x > \lambda^{1-\alpha}$ or $\alpha > \frac12$. 
\end{enumerate} 

This is essentially the construction  we shall employ below, up to an important
detail: We want to obtain an approximate solution in the energy space, which 
forces us to do technical modifications so that our functions are $x$ derivatives of
 suitable functions.

\section{Construction of the approximate solution}

We begin the construction  by collecting several elementary technical 
observations  needed to obtain good antiderivatives with respect to $x$.
If $f\in C_{0}^{\infty}(\R)$ is such that 
\begin{equation}\label{suq}
\int_{-\infty}^{\infty}f(x)dx=0,
\end{equation}
then, for every $x\in\R$,
$$
\big|(\partial_{x}^{-1}f)(x)\big|=
\big|\int_{-\infty}^{x}f(y)dy\big|
\leq {\rm mes}\big({\rm supp}(f)\big)\big(\sup_{y\in\R}|f(y)|\big).
$$
In particular, if for some $R>0$, ${\rm supp}(f)\subset [-R,R]$ then  for every
$x\in\R$,
$$
\big|
(\partial_{x}^{-1}f)(x)\big|\leq 2R\big(\sup_{y\in\R}|f(y)|\big),
$$
and if in addition
$$
\int_{-\infty}^{\infty}xf(x)dx=0,
$$
then we also have
$$
\big|
(\partial_{x}^{-2}f)(x)\big|\leq 4R^{2}\sup_{y\in\R}|f(y)|.
$$
Let us also notice that if $f\in C_{0}^{\infty}(\R)$ is such that (\ref{suq}) holds
and 
for some $R>0$, ${\rm supp}(f)\subset [-R,R]$ then
$$
{\rm supp}\big(\partial_{x}^{-1}(f)\big)\subset [-R,R].
$$
We recall 
\eqref{alphabeta}, 
that $\varphi \in C^\infty$ is supported in $[-2,2]$, identically $1$ 
in $[-1,1]$,   $\lambda>1$ and set 
$$
\psi_{\lambda}(x)=\varphi\big(\frac{x}{\lambda^{\alpha}}\big)
-2\varphi\big(\frac{x}{\lambda^{\alpha}}+c_{\lambda}\big)
+
\varphi\big(\frac{x}{\lambda^{\alpha}}+2c_{\lambda}\big),
$$
where
$$
c_{\lambda}=\frac{2\pi[10\lambda^{1+\alpha}]}{\lambda^{1+\alpha}}\,,
$$
with $[s]$ denoting the largest integer  $\le s$.
Notice that, $\psi_{\lambda}$ is supported in an interval of
size $\sim\lambda^{\alpha}$. 
In addition, for every $\gamma\in\R$,
$$
\int_{-\infty}^{\infty}\psi_{\lambda}(x)\cos(\lambda x+\gamma)dx=
\int_{-\infty}^{\infty}x\,\psi_{\lambda}(x)\cos(\lambda x+\gamma)dx=0\,.
$$
Therefore $\partial_{x}^{-1}(\psi_{\lambda}(x)\cos(\lambda x+\gamma))$ and 
$\partial_{x}^{-2}(\psi_{\lambda}(x)\cos(\lambda x+\gamma))$
are well defined $C_{0}^{\infty}(\R)$ functions.

Next, to shorten the notation, we define for $|\omega| \le 1 $
\begin{equation} 
\Phi_{\lambda}=\Phi_{\lambda}(t,x,y,\omega)=4\lambda^{3}t+\lambda
x+\sqrt{3}\lambda^{2}y+\omega t,
\end{equation} 
where we suppress $\omega$ in the notation of $\Phi_\lambda$, and we set 
\[  \widetilde \psi_\lambda(x) =
\varphi\big(\frac{x}{2\lambda^{\alpha}}\big)
-2\varphi\big(\frac{x}{2\lambda^{\alpha}}+c_{\lambda}/2\big)
+\varphi\big(\frac{x}{2\lambda^{\alpha}}+c_{\lambda}\big)\,
\] 
and 
\[ 
\varphi_\lambda(y)= \varphi(y/\lambda^\beta), 
\quad \widetilde \varphi_{\lambda}(y) = \varphi(y/(2\lambda^\beta)).
\]
For $|\omega|\leq 1$ and $\lambda\geq 1$, we define an approximate solution 
$u_{ap}$ of (\ref{1}) by the formula
\begin{equation}
\label{approximate}  
\begin{split} 
u_{ap}(t,x,y) = &
-\lambda^{-1-\frac{\alpha+\beta}{2}}\psi_{\lambda}(x)\varphi_{\lambda}(y)\cos(4\lambda^{3}t+\lambda
x+\sqrt{3}\lambda^{2}y+\omega t)
\\
 &
-\lambda^{-1}\omega\,\widetilde{\psi}_{\lambda}(x)\widetilde{\varphi}_{\lambda}(y)\,.
\end{split} 
\end{equation} 
Notice that
$$
\int_{-\infty}^{\infty}\widetilde{\psi}_{\lambda}(x)dx=\int_{-\infty}^{\infty}x\,\widetilde{\psi}_{\lambda}(x)dx=0\,.
$$
Therefore $\partial_{x}^{-1}(\widetilde{\psi}_{\lambda})$ and 
$\partial_{x}^{-2}(\widetilde{\psi}_{\lambda})$
are well defined $C_{0}^{\infty}(\R)$ functions.
Moreover, for $\lambda\gg 1$,
$$
\psi_{\lambda}\widetilde{\psi}_{\lambda}=\psi_{\lambda}.
$$
The main properties of $u_{ap}$ are collected in the following lemma.
\begin{lemme}\label{L2}
There exist $\delta>0$, $c>0$ and $C>0$ such that for every $\omega\in [-1,1]$,
every $\lambda\geq
1$,
\begin{equation}\label{I}
\Big\|
(\partial_{t}+\partial_{x}^{3}-\partial_{x}^{-1}\partial_{y}^{2})u_{ap}+u_{ap}\partial_{x}(u_{ap})
\Big\|_{L^2(\R^2)}\leq C\lambda^{-1-\delta}\,.
\end{equation}

Moreover
\begin{equation}\label{II}
\|\partial_{x}^{-1}\partial_{y}u_{ap}(t)\|_{L^2(\R^2)}\leq C,
\end{equation}
\begin{equation}\label{III}
\|\partial_{x}^{-2}\partial_{y}^{2}u_{ap}(t)\|_{L^2(\R^2)}\leq C\lambda
\end{equation}
and, for every $t,\omega, \omega^\prime\in [-1,1]$,
 \begin{equation}\label{IV} 
 \Vert  \partial_x (u_{ap,\omega}(t) -  u_{ap,\omega^\prime}(t)) 
\Vert_{L^2(\R^2)} \ge  c |\omega -\omega^\prime|| t|-C\lambda^{-\delta} .  
\end{equation} 
\end{lemme}
\begin{remark} It is not hard to keep track of the size of $\delta$. Let
$\varepsilon$ be a small positive constant, choose $\beta = 2 \alpha = \frac43 -
\varepsilon$. Then $\delta$ may be chosen to be $\frac13 - \varepsilon$.
\end{remark}

\begin{proof}
In the proof of this lemma, we denote by $o_{L^2} (\lambda^{-1})$ quantities having
$L^2(\R^2)$ norm bounded by $C\lambda^{-1-\delta}$ for a suitable $\delta>0$
uniformly for 
 $t\in[-1,1]$, $\omega\in[-1,1]$ and what is the most important, $\lambda\geq 1$.
The proof requires elementary but careful calculations. 

It is easy to check, using integration by parts,  that
$$
\|\partial_{x}^{-1}\partial_{y}^{2}\big(
\lambda^{-1}\omega\,\widetilde{\psi}_{\lambda}(x)\widetilde{\varphi}_{\lambda}(y)\big)\|_{L^2(\R^2)}\leq
C\lambda^{-1+\frac{\alpha+\beta}{2}+\alpha-2\beta}= C \lambda^{-1+ \frac32
(\alpha-\beta)}\,.
$$
Next,
$$
\|\partial_x^3\big(
\lambda^{-1}\omega\,\widetilde{\psi}_\lambda(x)\widetilde{\varphi}_{\lambda}(y)\big)\|_{L^2(\R^2)}\leq
C\lambda^{-1+\frac{\alpha+\beta}{2}-3\alpha}
$$
and thus, thanks to the assumptions on $(\alpha,\beta)$, we obtain that

\begin{equation}\label{2}
(\partial_{x}^{3}-\partial_{x}^{-1}\partial_{y}^{2})
\big(\lambda^{-1}\omega\,\widetilde{\psi}_{\lambda}(x)
\widetilde{\varphi}_{\lambda}(y)\big)=
O_{L^2} (\lambda^{-1+\frac{\alpha+\beta}2 + \max\{ \alpha-2\beta,-3\alpha\} })=
 o_{L^2}(\lambda^{-1})\,.
\end{equation}
Coming back to the definition of $u_{ap}$, we can readily check that
\begin{equation}\label{3}
u_{ap}\,\partial_{x}
u_{ap}=-\omega\lambda^{-1-\frac{\alpha+\beta}{2}}\psi_{\lambda}(x)\varphi_{\lambda}(y)\sin(\Phi_{\lambda}(t,x,y,\omega))
+O_{L^2}(\lambda^{-2-\frac{\alpha-\beta}2} ).
\end{equation} 

Notice that the leading term in (\ref{3}) is coming from the product of the
high frequency part of $\partial_{x}u_{ap}$ and the low frequency part of
$u_{ap}$.

Next, we compute integrating by parts
\begin{equation*}
\partial_{x}^{-1}(\psi_{\lambda}\cos\Phi_{\lambda})
=\int_{-\infty}^{x}\psi_{\lambda}\cos\Phi_{\lambda}=
\lambda^{-1}\psi_{\lambda}\sin\Phi_{\lambda}-\lambda^{-1}\int_{-\infty}^{x}\partial_{x}[\psi_{\lambda}]\sin\Phi_{\lambda}\,.
\end{equation*}
We integrate by parts two more times to arrive at
\begin{multline}\label{4}
\partial_{x}^{-1}(\psi_{\lambda}\cos\Phi_{\lambda})=
\lambda^{-1}\psi_{\lambda}\sin\Phi_{\lambda}
+
\lambda^{-2}\partial_{x}[\psi_{\lambda}]\cos\Phi_{\lambda}-
\\
-
\lambda^{-3}\partial_{x}^{2}[\psi_{\lambda}]\sin\Phi_{\lambda}
+
\lambda^{-3}\int_{-\infty}^{x}
\partial_{x}^{3}[\psi_{\lambda}]\sin\Phi_{\lambda}.
\end{multline}
Using the Leibniz rule, since $\beta>1$, we infer that
$$
\lambda^{-1-\frac{\alpha+\beta}{2}}
\partial_{y}^{2}
\big([\lambda^{-1}\psi_{\lambda}\sin\Phi_{\lambda}]
\varphi_{\lambda}(y)\big)
=
-3\lambda^{2-\frac{\alpha+\beta}{2}}\psi_{\lambda}(x)\varphi_{\lambda}(y)\sin\Phi_{\lambda}+O_{L^2}(\lambda^{-\beta})\,.
$$
Similarly
\[
\begin{split} 
\lambda^{-1-\frac{\alpha+\beta}{2}}
\partial_{y}^{2}
\big(
[\lambda^{-2}\partial_{x}[\psi_{\lambda}]\cos\Phi_{\lambda}]
\varphi_{\lambda}(y)\big)= & 
-3\lambda^{1-\frac{\alpha+\beta}{2}}\partial_{x}[\psi_{\lambda}(x)]\varphi_{\lambda}(y)\cos\Phi_{\lambda}
\\ & +O_{L^2}(\lambda^{-1-\alpha-\beta})
\end{split} 
\]
and, 
$$
\lambda^{-1-\frac{\alpha+\beta}{2}}
\partial_{y}^{2}
\big(
[\lambda^{-3}\partial_{x}^{2}[\psi_{\lambda}]\sin\Phi_{\lambda}]
\varphi_{\lambda}(y)\big)=O_{L^2}(\lambda^{-2\alpha}).
$$
We recall that $\alpha> \frac12$.
Similarly
$$
\lambda^{-4 -\frac{\alpha+\beta}2} \partial_y^2 \int_{-\infty}^{x}
\partial_{x}^{3}[\psi_{\lambda}]\sin\Phi_{\lambda} = O_{L^2}( \lambda^{-2\alpha}).
$$
Summarizing, we can conclude that
\begin{equation}
\begin{split}
\partial_{x}^{-1}\partial_{y}^{2}
\big(
\lambda^{-1-\frac{\alpha+\beta}{2}}\psi_{\lambda}(x)\varphi_{\lambda}(y)\cos\Phi_{\lambda}
\big)
 = &
-3\lambda^{2-\frac{\alpha+\beta}{2}}\psi_{\lambda}(x)\varphi_{\lambda}(y)\sin\Phi_{\lambda}
\\
 &
-3\lambda^{1-\frac{\alpha+\beta}{2}}\partial_{x}[\psi_{\lambda}(x)]\varphi_{\lambda}(y)\cos\Phi_{\lambda}
+o_{L^2}(\lambda^{-1}).
\end{split}
\end{equation}
Using the Leibniz rule, we infer 
\begin{eqnarray*}
\partial_{x}^{3}
\big(
\lambda^{-1-\frac{\alpha+\beta}{2}}\psi_{\lambda}(x)\varphi_{\lambda}(y)\cos\Phi_{\lambda}
\big)
& = &
\lambda^{2-\frac{\alpha+\beta}{2}}\psi_{\lambda}(x)\varphi_{\lambda}(y)\sin\Phi_{\lambda}
\\
& &
-3\lambda^{1-\frac{\alpha+\beta}{2}}\partial_{x}[\psi_{\lambda}(x)]\varphi_{\lambda}(y)\cos\Phi_{\lambda}
+o_{L^2}(\lambda^{-1}).
\end{eqnarray*}
Therefore using (\ref{2}), we obtain that
$$
(\partial_{x}^{3}-\partial_{x}^{-1}\partial_{y}^{2})u_{ap}=
-4\lambda^{2-\frac{\alpha+\beta}{2}}\psi_{\lambda}(x)\varphi_{\lambda}(y)\sin\Phi_{\lambda}
+o_{L^2}(\lambda^{-1}).
$$
The cancellation of the term
$$
3\lambda^{1-\frac{\alpha+\beta}{2}}\partial_{x}[\psi_{\lambda}(x)]\varphi_{\lambda}(y)\cos\Phi_{\lambda}
$$
is the  main new point in this paper. It is an analytic expression of the fact that
the $x$ component of the velocity vector vanishes for the plane wave which we have
chosen. Here, we essentially use that we are dealing
with the KP-I equation, i.e. the sign in front of
$\partial_{x}^{-1}\partial_{y}^{2}$
is crucial to achieve this cancellation.
Since
$$
\partial_{t}u_{ap}=4\lambda^{2-\frac{\alpha+\beta}{2}}\psi_{\lambda}(x)\varphi_{\lambda}(y)\sin\Phi_{\lambda}
+\lambda^{-1-\frac{\alpha+\beta}{2}}\omega\psi_{\lambda}(x)\varphi_{\lambda}(y)\sin\Phi_{\lambda},
$$
we obtain that
$$
(\partial_{t}+\partial_{x}^{3}-\partial_{x}^{-1}\partial_{y}^{2})u_{ap}
=
\omega\lambda^{-1-\frac{\alpha+\beta}{2}}\psi_{\lambda}(x)\varphi_{\lambda}(y)\sin\Phi_{\lambda}+o_{L^2}(\lambda^{-1}).
$$
This, together with  (\ref{3}), completes the proof of (\ref{I}).

Using (\ref{4}), arguing in the same way as there, we can write
\begin{equation} 
\partial_{x}^{-1}\partial_{y}
\big(
\lambda^{-1-\frac{\alpha+\beta}{2}}\psi_{\lambda}(x)\varphi_{\lambda}(y)\cos\Phi_{\lambda}
\big)
= 
\sqrt{3}\lambda^{-\frac{\alpha+\beta}{2}}\psi_{\lambda}(x)\varphi_{\lambda}(y)\cos\Phi_{\lambda}
+o_{L^2}(\lambda^{-1}).
\end{equation} 
Moreover
\begin{equation}
\|\partial_{x}^{-1}\partial_{y}\big(
\lambda^{-1}\omega\,\widetilde{\psi}_{\lambda}(x)\widetilde{\varphi}_{\lambda}(y)\big)\|_{L^2(\R^2)}\leq
C\lambda^{-1+\frac{\alpha+\beta}{2}+\alpha-\beta}\leq C
\end{equation}
which completes the proof of (\ref{II}).

Let us now turn to the proof of (\ref{III}). The low frequency part of $u_{ap}$
can be estimated as
$$
\|\partial_{x}^{-2}\partial_{y}^{2}\big(
\lambda^{-1}\omega\,\widetilde{\psi}_{\lambda}(x)\widetilde{\varphi}_{\lambda}(y)\big)\|_{L^2(\R^2)}\leq
C\lambda^{-1+\frac{\alpha+\beta}{2}+2\alpha-2\beta}\leq C
$$
We next estimate the high frequencies and repeat the calculation  
 of (\ref{4}),
\begin{multline*}
\partial_{x}^{-1}(\psi_{\lambda}\sin\Phi_{\lambda})
=\int_{-\infty}^{x}\psi_{\lambda}\sin\Phi_{\lambda}=
-\lambda^{-1}\psi_{\lambda}\cos\Phi_{\lambda}+\lambda^{-1}\int_{-\infty}^{x}\partial_{x}[\psi_{\lambda}]\cos\Phi_{\lambda}
\\
=
-\lambda^{-1}\psi_{\lambda}\cos\Phi_{\lambda}+
\lambda^{-2}\partial_{x}[\psi_{\lambda}]\sin\Phi_{\lambda}
-
\lambda^{-2}\int_{-\infty}^{x}\partial_{x}^{2}[\psi_{\lambda}]\sin\Phi_{\lambda}\,.
\end{multline*}
Next, we estimate each term in the right hand-side of the above equality. First
\begin{equation}\label{mar1}
\|\lambda^{-2}\psi_{\lambda}\cos\Phi_{\lambda}\|_{L^2(\R_{x})}
\leq C\lambda^{-2}\lambda^{\frac{\alpha}{2}}
\end{equation}
and then
\begin{equation}\label{mar2}
\|\lambda^{-3}\partial_{x}[\psi_{\lambda}]\sin\Phi_{\lambda}\|_{L^2(\R_{x})}
\leq C\lambda^{-3}\lambda^{-\alpha}\lambda^{\frac{\alpha}{2}}
\leq C\lambda^{-2}\lambda^{\frac{\alpha}{2}}
\end{equation}
and finally
$$
\|\lambda^{-3}\int_{-\infty}^{x}\partial_{x}^{2}[\psi_{\lambda}]\sin\Phi_{\lambda}\|_{L^2(\R_{x})}
\leq C\lambda^{-3}\lambda^{-2\alpha}\lambda^{\alpha}\lambda^{\frac{1}{2}}
\leq C\lambda^{-2}\lambda^{\frac{\alpha}{2}}.
$$
Notice that (\ref{mar1}) and (\ref{mar2}) imply that
$$
\|\partial_{x}^{-1}(\lambda^{-2}\partial_{x}[\psi_{\lambda}]\cos\Phi_{\lambda})\|_{L^2(\R_{x})}\leq
C\lambda^{-2}\lambda^{\frac{\alpha}{2}}
$$
which is the relevant bound for the second term in the right hand-side of (\ref{4}).
It remains to estimate the last two terms in the right hand-side of (\ref{4}).
We can write
$$
\|\partial_{x}^{-1}[\lambda^{-3}\partial_{x}^{2}[\psi_{\lambda}]\sin\Phi_{\lambda}]\|_{L^2(\R_{x})}
\leq C\lambda^{-3}\lambda^{-2\alpha}
\lambda^{\alpha}\lambda^{\frac{1}{2}}\leq
C\lambda^{-2}\lambda^{\frac{\alpha}{2}},
$$
since $\alpha>1/2$. 
For the last term in the right hand-side of (\ref{4}), we can write
\begin{eqnarray*}
\big\|
\partial_{x}^{-1}[\lambda^{-3}\int_{-\infty}^{x}\partial_{x}^{3}[\psi_{\lambda}]\sin\Phi_{\lambda}]
\big\|_{L^2(\R_{x})}
& = &
\big\|
\partial_{x}^{-2}[\lambda^{-3}\partial_{x}^{3}[\psi_{\lambda}]\sin\Phi_{\lambda}]
\big\|_{L^2(\R_{x})}
\\
& \leq & C\lambda^{-3}\lambda^{-3\alpha}\lambda^{\frac{1}{2}}\lambda^{2}
\leq C\lambda^{-2}\lambda^{\frac{\alpha}{2}}
\end{eqnarray*}
since $\alpha>1/2>3/7$. 

Summarizing, we infer that the high frequencies of $u_{ap}$ can be estimated as 
$$
\|\partial_{x}^{-2}(\psi_{\lambda}\cos\Phi_{\lambda})\|_{L^2(\R_{x})}
\leq C\lambda^{-2}\lambda^{\frac{\alpha}{2}}
$$
and thus, using that $\partial_{y}^{2}$ is causing at most an amplification
factor $\lambda^4$, we conclude that
$$
\big\|\partial_{x}^{-2}\partial_{y}^{2}
\big(
\lambda^{-1-\frac{\alpha+\beta}{2}}\psi_{\lambda}(x)\varphi_{\lambda}(y)\cos\Phi_{\lambda}
\big)
\big\|_{L^2(\R^2)}\leq
C\lambda^{-1-\frac{\alpha+\beta}{2}}(\lambda^{-2}\lambda^{\frac{\alpha}{2}})
(\lambda^{4}\lambda^{\frac{\beta}{2}})
=C\lambda.
$$
This proves (\ref{III}).

Finally, we give the proof of \eqref{IV}. Notice that
\begin{equation}\label{dx}  
\partial_x u_{ap} = \lambda^{-\frac{\alpha+\beta}2}\psi_\lambda(x)
\varphi_\lambda(y) \sin( \Phi_\lambda)   + O_{L^2} (
\lambda^{-1+\frac{\beta-\alpha}2}).
\end{equation}
With $a =  4\lambda^3 t + \lambda x + \sqrt3 \lambda^2 y$, we may write
\[ 
\sin (a + \omega t)
-  \sin (a+ \omega^\prime t)
 = 
2\sin(t(\omega-\omega')/2)
\cos\big(a+t(\omega+\omega')/2), 
\]
and after a sequence of integrations by parts, we get
\[ 
\Vert \lambda^{-\frac{\alpha+\beta}2} \psi_\lambda \varphi_\lambda 
\big\{ \sin \Phi_{\omega,\lambda}  
- \sin \Phi_{\omega,\lambda^\prime}  \big\} 
\Vert_{L^2}^2 
\geq c(|t||\omega-\omega'|)^{2}  
\Vert 
\lambda^{-\frac{\alpha+\beta}2} \psi_\lambda \varphi_\lambda \Vert_{L^2} ^{2}
-C\lambda^{-2} . 
\] 
Using the choice of $c_{\lambda}$, we can minorate
$\|\psi_{\lambda}\|_{L^2(\R)}$ and thus
\begin{equation*} 
\Vert \lambda^{-\frac{\alpha+\beta}2} \psi_\lambda \varphi_\lambda 
  \big\{ \sin \Phi_{\omega,\lambda}  - \sin \Phi_{\omega^\prime, \lambda}\big\}
\Vert_{L^2} 
  \ge c |\omega-\omega^\prime| |t| 
-C\lambda^{-1} 
\end{equation*} 
which proves (\ref{IV}).
This completes the proof of Lemma~\ref{L2}.
\end{proof}

\section{Bounds for the exact solution} 
Let $u_{\omega,\lambda}(t,x,y)$ be a solution of the KP-I equation with data
 \begin{equation*}
u_{\omega,\lambda}(0,x,y)=-\lambda^{-1-\frac{\alpha+\beta}{2}}\psi_{\lambda}(x)\varphi_{\lambda}(y)\cos(\lambda
x+\sqrt{3}\lambda^{2}y)
-\lambda^{-1}\omega\,\widetilde{\psi}_{\lambda}(x)\widetilde{\varphi}_{\lambda}(y).
\end{equation*}
Thanks to the properties of $\psi_{\lambda}(x)$ and
$\widetilde{\psi}_{\lambda}(x)$, we can apply the global well-posedness result
of \cite{K} to obtain that $u_{\omega,\lambda}(t,x,y)$ is globally defined and
satisfies the conservation laws mentioned in the introduction (see also
\cite[Proposition~4]{MST1}).
Moreover, for every $t$,
$\xi^{-1}\widehat{u_{\omega,\lambda}}(t,\xi,\eta)$ belongs to $L^2(\R^2)$
(see \cite{IN}).

In order to bound $u_{\omega,\lambda}$ in $Z$, we will use the following 
anisotropic Sobolev inequality.
\begin{lemme}\label{sobolev}
For $2\leq p\leq 6$ there exists $C>0$ such that for every $u\in X$,
\begin{equation}\label{malak}
\|u\|_{L^p(\R^2)}\leq C\|u\|_{L^2(\R^2)}^{\frac{6-p}{2p}}\,\,
\|u_x\|_{L^2(\R^2)}^{\frac{p-2}{p}}\,\,
\|\partial_x^{-1}u_y\|_{L^2(\R^2)}^{\frac{p-2}{2p}}\,\, .
\end{equation}
\end{lemme}
We refer to \cite{Tom} for a proof of (\ref{malak}).
The $L^2$ conservation law yields,
$$
\|u_{\omega,\lambda}(t,\cdot)\|_{L^2(\R^2)}
=
\|u_{\omega,\lambda}(0,\cdot)\|_{L^2(\R^2)}
=
\|u_{ap}(0,\cdot)\|_{L^2(\R^2)}
\leq C\,.
$$
The energy conservation, (\ref{II}) and (\ref{malak}) with $p=3$ yield
$$
E(u_{\omega,\lambda}(t,\cdot))=E(u_{\omega,\lambda}(0,\cdot))=
E(u_{ap}(0,\cdot))\leq C\,.
$$
Another use of (\ref{malak}) with $p=3$ then gives the following bound for
the leading part of the energy,
\begin{equation}\label{energy}
\|\partial_{x}(u_{\omega,\lambda})(t,\cdot)\|_{L^2(\R^2)}+\|\partial_{x}^{-1}\partial_{y}(u_{\omega,\lambda})(t,\cdot)\|_{L^2(\R^2)}\leq
C.
\end{equation}
We now establish several bounds for the cubic and quartic terms of the
functional $F$ of \eqref{F}. We can write, by invoking (\ref{malak}) with $p=3,4$
\begin{eqnarray*}
\Big|\frac{5}{6} \int_{\R^2} u^2 (\partial^{-2}_x u_{yy})\Big|
& \leq & \|u\|_{L^4(\R^2)}^{2}\|\partial^{-2}_x u_{yy}\|_{L^2(\R^2)}
\\
& \leq  &
C\|u\|_{L^2(\R^2)}^{\frac{1}{2}}
\|u_x\|_{L^2(\R^2)}
\|\partial^{-1}_x u_{y}\|_{L^2(\R^2)}^{\frac{1}{2}}
\|\partial^{-2}_x u_{yy}\|_{L^2(\R^2)}
\end{eqnarray*}
and
\[
\begin{split} 
\Big|\int_{\R^2} u\, (\partial_{x}^{-1} u_y)^2\Big|
& \leq 
\|u\|_{L^3(\R^2)}
\|\partial^{-1}_x u_{y}\|_{L^3(\R^2)}^{2}
\\
& \leq 
C\|u\|_{L^2(\R^2)}^{\frac{1}{2}}
\|u_x\|_{L^2(\R^2)}^{\frac{1}{3}}
\|\partial^{-1}_x u_{y}\|_{L^2(\R^2)}^{\frac{7}{6}}
\|u_y\|_{L^2(\R^2)}^{\frac{2}{3}}
\|\partial^{-2}_x u_{yy}\|_{L^2(\R^2)}^{\frac{1}{3}}\,.
\end{split} 
\]
Next,
\begin{eqnarray*}
\Big|\int_{\R^2} u^2\, u_{xx}\Big| 
& \leq &
\|u\|_{L^4(\R^2)}^{2}
\|u_{xx}\|_{L^2(\R^2)}
\\
& \leq &
C\|u\|_{L^2(\R^2)}^{\frac{1}{2}}\|u_x\|_{L^2(\R^2)}
\|\partial^{-1}_x u_{y}\|_{L^2(\R^2)}^{\frac{1}{2}}
\|u_{xx}\|_{L^2(\R^2)}
\end{eqnarray*}
and finally
$$
\Big|\int_{\R^2} u^4\Big|\leq C
\|u\|_{L^2(\R^2)}\|u_x\|_{L^2(\R^2)}^{2}
\|\partial^{-1}_x u_{y}\|_{L^2(\R^2)}\,.
$$
Using the above bounds, estimates (\ref{III}), (\ref{energy}), and the
conservation of $F$, we obtain that for $\lambda\geq 1$,
$$
F(u_{\omega,\lambda}(t,\cdot))=F(u_{\omega,\lambda}(0,\cdot))=
F(u_{ap}(0,\cdot))\leq C\lambda^2\,.
$$
Using again the estimates for the cubic and the quartic terms of $F$, we
obtain that the leading part of $F$ satisfies for $\lambda\geq 1$,
$t\in[-1,1]$,
$$
\|\partial_{x}^{2}(u_{\omega,\lambda})(t,\cdot)\|_{L^2(\R^2)}+
\|\partial_{y}(u_{\omega,\lambda})(t,\cdot)\|_{L^2(\R^2)}+
\|\partial_{x}^{-2}\partial_{y}^{2}(u_{\omega,\lambda})(t,\cdot)\|_{L^2(\R^2)}\leq
C\lambda.
$$
\section{The difference between  approximate and exact solution}

We begin by controlling the size of $u_{ap}$.
Using (\ref{III}) and the definition of $u_{ap}$, we infer that for
$\lambda\geq 1$, $t\in[-1,1]$,
$$
\|\partial_{x}^{2}(u_{ap})(t,\cdot)\|_{L^2(\R^2)}+
\|\partial_{y}(u_{ap})(t,\cdot)\|_{L^2(\R^2)}+
\|\partial_{x}^{-2}\partial_{y}^{2}(u_{ap})(t,\cdot)\|_{L^2(\R^2)}\leq C\lambda
$$
and thus, with 
\begin{equation}
v_{\omega,\lambda}= u_{\omega,\lambda}-u_{ap},
\end{equation}
$$
\|\partial_{x}^{2}(v_{\omega,\lambda})(t,\cdot)\|_{L^2(\R^2)}+
\|\partial_{y}(v_{\omega,\lambda})(t,\cdot)\|_{L^2(\R^2)}+
\|\partial_{x}^{-2}\partial_{y}^{2}(v_{\omega,\lambda})(t,\cdot)\|_{L^2(\R^2)}\leq
C\lambda.
$$
In particular
\begin{equation}\label{dve}
\|\partial_{x}^{2}(v_{\omega,\lambda})(t,\cdot)\|_{L^2(\R^2)}\leq C\lambda.
\end{equation}
We next bound the $L^2$ norm of $v_{\omega,\lambda}$.
\begin{lemme}\label{nula}
There exist $\delta>0$  such that 
\begin{equation}\label{zero}
\|v_{\omega,\lambda}(t,\cdot)\|_{L^2(\R^2)}\leq C\lambda^{-1-\delta}.
\end{equation}
uniformly in $\lambda\ge 1$, $|\omega|\le 1$ and $|t|\le 1$.
\end{lemme}
\begin{proof}
The function $v_{\omega,\lambda}$ solves the equation
\begin{equation}\label{star}
(\partial_{t}+\partial_{x}^{3}-\partial_{x}^{-1}\partial_{y}^{2})v_{\omega,\lambda}+
v_{\omega,\lambda}\partial_{x}(v_{\omega,\lambda})+\partial_{x}(u_{ap}v_{\omega,\lambda})+G=0,
\end{equation}
where $v_{\omega,\lambda}(0,x,y)=0$ and
$$
G=(\partial_{t}+\partial_{x}^{3}-\partial_{x}^{-1}\partial_{y}^{2})u_{ap}+u_{ap}\partial_{x}(u_{ap})\,.
$$
Thanks to (\ref{I}),
$$
\|G(t,\cdot)\|_{L^2(\R^2)}\leq C\lambda^{-1-\delta},\quad \delta>0.
$$
Multiplying (\ref{star}) by $v_{\omega,\lambda}$ and an integration over
$\R^2$ gives
\begin{eqnarray*}
\frac{d}{dt}\|v_{\omega,\lambda}(t,\cdot)\|_{L^2(\R^2)}^{2}
& \lesssim &
\|\partial_x
u_{ap}(t,\cdot)\|_{L^{\infty}(\R^2)}\|v_{\omega,\lambda}(t,\cdot)\|_{L^2(\R^2)}^{2}
\\
& &
+
\|v_{\omega,\lambda}(t,\cdot)\|_{L^2(\R^2)}\|G(t,\cdot)\|_{L^2(\R^2)}\,.
\end{eqnarray*}
>From the definition of $u_{ap}$, we infer that
$$
\|\partial_x u_{ap}(t,\cdot)\|_{L^{\infty}(\R^2)}\leq C\lambda^{-1}\,.
$$
Therefore, by Gronwall's inequality
for $t\in[-1,1]$,
$$
\|v_{\omega,\lambda}(t,\cdot)\|_{L^2(\R^2)}\leq
\sup_{t\in[-1,1]}\|F(t,\cdot)\|_{L^2(\R^2)}\leq C\lambda^{-1-\delta}\,.
$$
This completes the proof of Lemma~\ref{nula}.
\end{proof}
Interpolation between (\ref{dve}) and (\ref{zero}) gives that for 
$\lambda\geq 1$,
\begin{equation}\label{edno}
\|\partial_{x}(v_{\omega,\lambda})(t,\cdot)\|_{L^2(\R^2)}\leq C\lambda^{-\delta/2}.
\end{equation}
After these preparation we turn to the proof of Theorem~\ref{th1}. 
Consider the two families of solutions $(u_{1,\lambda})$ and
$(u_{-1,\lambda})$, $\lambda\gg 1$. Write for $\lambda\geq 1$,
\begin{eqnarray*}
\|u_{1,\lambda}(0,\cdot)-u_{-1,\lambda}(0,\cdot)\|_{X} & = &2
\|\lambda^{-1}\omega\,\widetilde{\psi}_{\lambda}(x)\widetilde{\varphi}_{\lambda}(y)\|_{X}
\\
& \le &
2\lambda^{-1}\|\widetilde{\psi}_{\lambda}(x)\widetilde{\varphi}_{\lambda}(y)\|_{L^2(\R^2)}
\\
& &
+2\lambda^{-1}\|\partial_{x}[\widetilde{\psi}_{\lambda}(x)]\widetilde{\varphi}_{\lambda}(y)\|_{L^2(\R^2)}
\\
& &
+2\lambda^{-1}\|\partial_{x}^{-1}[\widetilde{\psi}_{\lambda}(x)]\partial_{y}[\widetilde{\varphi}_{\lambda}(y)]\|_{L^2(\R^2)}
\\
&\leq&
C\lambda^{-1}\lambda^{\frac{\alpha+\beta}{2}}+C\lambda^{-1}\lambda^{\alpha}\lambda^{-\beta}\lambda^{\frac{\alpha+\beta}{2}}\,.
\end{eqnarray*}
Thanks to the assumptions on $(\alpha,\beta)$, we obtain that
$$
\lim_{\lambda\rightarrow\infty}\|u_{1,\lambda}(0,\cdot)-u_{-1,\lambda}(0,\cdot)\|_{X}=0\,.
$$
To conclude we  provide a nontrivial lower bound on 
\[ \liminf_{\lambda \to \infty} \Vert \partial_x ( u_{1,\lambda}-u_{-1,\lambda})
\Vert_{L^2}.\]
Equation \eqref{edno} reduces this to the corresponding statement for $u_{ap}$,
which is inequality \eqref{IV}.
This completes the proof of Theorem~\ref{th1}.
\qed
\begin{remarque}
Actually, we proved a stronger statement than Theorem~\ref{th1}.
We obtained the existence of two families of solutions of the KP-I equation
which remain bounded in the energy space, such that their difference tend to
zero in the energy space but such that for $t\in[-1,1]$, $t\neq 0$ the $x$
derivative of their difference in $L^2(\R^2)$, which is only a part of the
energy norm, does not tend to zero.
\end{remarque}
\begin{remarque}
As in \cite{KT}, if one is interested to show the failure of uniform
continuity of the flow of KP-I on $H^s(\R^2)$ for large $s$ a modification of
the low frequency part of the approximate solution is needed. Namely one should replace
\begin{equation}\label{lowpart}
\omega\lambda^{-1}\widetilde{\psi}_{\lambda}(x)
\widetilde{\varphi}_{\lambda}(y)
\end{equation}
in $u_{ap}$ by the solution of the KP-I equation with initial data (\ref{lowpart}).
We refer to \cite{KT} for the details of this construction.
\end{remarque}

\end{document}